\documentclass[12pt]{article}
\usepackage{amssymb}
\usepackage{amsthm}
\usepackage{amsmath}
\usepackage{graphicx}
 \newtheorem{thm}{Theorem}[section]

 \theoremstyle{definition}
 
 \newtheorem{rem}[thm]{Remark}
 \numberwithin{equation}{section}

\begin{document}
\title{Global Mild Solutions of the Navier-Stokes
Equations \footnote{Corresponding authors:
Z. Lei, School of Mathematical Sciences, Fudan Unviersity; F.-H.
Lin, Courant Institute of Mathematics, New
  York University}}
\author{Zhen Lei\footnote{School of Mathematical Sciences; LMNS and Shanghai Key
  Laboratory for Contemporary Applied Mathematics, Fudan
  University, Shanghai 200433, P. R. China. {\it Email:
  leizhn@yahoo.com}}
\and Fang-Hua Lin\footnote{Courant Institute of Mathematics, New
  York University, USA. {\it Email: linf@cims.nyu.edu}}}
\date{\today}
\maketitle

\begin{abstract}
Here we establish a global well-posedness of \textit{mild} solutions
to the three-dimensional incompressible Navier-Stokes equations if the
initial data are in the space $\mathcal{X}^{-1}$ defined by $(1.3)$
and if the norms of the initial
data in $\mathcal{X}^{-1}$ are bounded exactly by the viscosity
coefficient $\mu$.
\end{abstract}


\section{The Result}


The incompressible Navier-Stokes equations in $\mathbb{R}_+ \times
\mathbb{R}^3$ are:
\begin{equation}\label{NS}
\begin{cases}
\partial_tv + v\cdot\nabla v + \nabla p = \mu\Delta v,\quad t > 0, x \in \mathbb{R}^3,\\[-3mm]\\
\nabla\cdot v = 0,\quad t > 0, x \in \mathbb{R}^3,
\end{cases}
\end{equation}
where $v$ is the velocity field of the fluid, $p$ is the pressure
and the constant $\mu$ is the viscosity.
To solve the Navier-Stokes
equations \eqref{NS} in $\mathbb{R}_+ \times \mathbb{R}^3$, one
assumes that the initial data
\begin{equation}\label{Data}
v(0, x) = v_0(x)
\end{equation}
are divergent free and possess certain regularity.


The global existence of weak solutions goes back to Leray \cite{Leray34}
and Hopf \cite{Hopf51}.
The global well-posedness of strong solutions for small initial
data is due to Fujita and Kato \cite{FK} (see also Chemin
\cite{Chemin}) in the Sobolev spaces ${\rm \dot{H}}^s$, $s \geq
\frac{1}{2}$, Kato \cite{Kato} in the Lebesgue space
$L^3(\mathbb{R}^3)$,
and Koch and Tataru \cite{KT} in the space ${\rm BMO}^{- 1}$ (see
also \cite{Cannone} and \cite{Planchon}). It should be noted, in
all these works, that the norms in corresponding spaces of the
initial data are assumed to be very small, say they are smaller
than the viscosity coefficient $\mu$ multiplied by a tiny positive
constant $\epsilon$.

In this note, we shall prove a new
global well-posedness result
for the three-dimensional incompressible Navier-Stokes equations.
Our \textit{mild} solutions will be in a scale-invariant function space
which is natural with respect to the scalings of the Navier-Stokes
equations. Moreover, we show that it is sufficient to assume the norms
of the initial data to be less than the  viscosity coefficient $\mu$.
The function space we will use here is
\begin{equation}\label{1.3}
\mathcal{X}^{-1} = \{f \in \mathcal{D}^\prime(R^3):
\int_{\mathbb{R}^3}|\xi|^{-1}|\widehat{f}|d\xi < \infty\}.
\end{equation}
Here $\mathcal{D}^\prime(R^3)$ represents the space of
distributions and $\widehat{f}$ represents the Fourier transform
of $f$. The norm of $\mathcal{X}^{-1}$ will be denoted by
$\|\cdot\|_{\mathcal{X}^{-1}}$. We will also use the notation
\begin{equation}\nonumber
\mathcal{X}^1 = \{f \in \mathcal{D}^\prime(R^3):
\int_{\mathbb{R}^3}|\xi||\widehat{f}|d\xi < \infty\}.
\end{equation}
The norm of $\mathcal{X}^1$ is denoted by
$\|\cdot\|_{\mathcal{X}^1}$.

\begin{thm}\label{thm1}
The three-dimensional incompressible Navier-Stokes equation
\eqref{NS} is well-posed globally in time for the initial data in
$\mathcal{X}^{-1}$ satisfying
\begin{equation}\label{1.4}
\|v_0\|_{\mathcal{X}^{-1}}  < \mu.
\end{equation}
Moreover, the solution $v$ is in $C(\mathbb{R}_+;
\mathcal{X}^{-1}) \cap L^1(\mathbb{R}_+; \mathcal{X}^1)$ and
satisfies
\begin{equation}\nonumber
\sup_{0 \leq t < \infty}\Big(\|v(t)\|_{\mathcal{X}^{-1}}  +
\big(\mu - \|v_0\|_{\mathcal{X}^{-1}}\big)\int_0^t\|\nabla
v(\tau)\|_{L^\infty}d\tau\Big) \leq \|v_0\|_{\mathcal{X}^{-1}}.
\end{equation}
\end{thm}


Let us remark that the space ${BMO^{- 1}}$ is the largest space
which is included in the tempered distribution and translation and
scaling invariant, (see, for instance, the paper by Chemin and
Gallagher \cite{CG}). Our space $\mathcal{X}^{-1}$ is contained in
${BMO^{- 1}}$. In fact, for each $f \in \mathcal{X}^{-1}$, write
$f = \nabla\cdot(\nabla\Delta^{-1}f)$. It is clear that
$\|\nabla\Delta^{-1}f\|_{{\rm BMO}} \leq
\|\nabla\Delta^{-1}f\|_{L^\infty} \leq
\int|\xi|^{-1}|\widehat{f}|d\xi = \|f\|_{\mathcal{X}^{-1}}$.
Hence, by Theorem 1 in \cite{KT}, one has $f \in {BMO^{- 1}}$. Let
us also mention the classical example of the initial data given in
\cite{CG}
\begin{equation}\nonumber
u_0^\epsilon(x) =
\frac{1}{\epsilon}\cos\frac{x_3}{\epsilon}(\partial_2\phi, -
\partial_1\phi, 0)^T.
\end{equation}
Due to Lemma 3.1 in \cite{CG}, it is easy to see that
$\|u_0^\epsilon\|_{\mathcal{X}^{-1}}$ is uniformly bounded
independent of $\epsilon$.

It is not hard to see that if the initial data are in the Sobolev
space $H^s$, $s > \frac{1}{2}$, then they are also in the space
$\mathcal{X}^{-1}$. However, this turns out to be false for $s =
\frac{1}{2}$. An easy counterexample is given as follows: Let $f$
be a non-negative function in the Schwartz class
$\mathcal{S}(\mathbb{R}^3)$ which is supported in the set $\{\xi
\in \mathbb{R}^3: 1 < |\xi| < 2\}$. Consider  $g$ which is defined
by
\begin{equation}\nonumber
\widehat{g} = \sum_{j \geq 1}\frac{2^{-2j}}{j}f(2^{-j}\xi).
\end{equation}
It is clear that $g \in \dot{H}^{\frac{1}{2}}$ but
\begin{equation}\nonumber
\|g\|_{\mathcal{X}^{-1}} = \sum_{j \geq
1}\int\frac{2^{-3j}}{j}f(2^{-j}\xi)d\xi = \|f\|_{L^1}\sum_{j \geq
1}\frac{1}{j} = \infty.
\end{equation}

\section{Proof of Theorem}

Let $\zeta$ be the standard mollifier in $\mathbb{R}^3$: $\zeta
\in C_0^\infty$, $0 \leq \zeta \leq 1$, $\int\zeta(x) dx = 1$. For
$\lambda > 0$, let $\zeta^\lambda(x) =
\lambda^{-3}\zeta(\lambda^{-1}x)$ and $v_0^\lambda = \zeta^\lambda
\ast v_0$. For $v_0 \in \mathcal{X}^{-1}$, since
$|\widehat{\zeta}(\xi)| \leq \int\zeta(x)dx = 1$, one has
\begin{equation}\label{2.1}
\begin{cases}
\|v_0^\lambda\|_{\mathcal{X}^{ -1}} = \int|\xi|^{-
1}|\widehat{v_0}(\xi)||\widehat{\zeta}(\lambda\xi)|d\xi \leq
\|v_0\|_{\mathcal{X}^{ -1}},\\[-4mm]\\
\|v_0^\lambda\|_{L^\infty} \lesssim
\int|\xi||\widehat{\zeta}(\lambda\xi)||\xi|^{-
1}|\widehat{v_0}(\xi)|d\xi \leq C_\lambda\|v_0\|_{\mathcal{X}^{
-1}}.
\end{cases}
\end{equation}
Consequently, by the standard local existence theory of the
Navier-Stokes equations, there exists a unique local smooth
solution $v^\lambda(t, x)$ on some time internal $[0, T_\lambda)$.
The associated pressure $p^\lambda$ is given by $p^\lambda =
(R\otimes R):(v\otimes v)$ with $R$ being the Riesz operator
(see, for instance, \cite{GIM}).

We need to derive an \textit{a priori} estimate under the
condition \eqref{1.4}. First of all, it is easy to see that
\begin{equation}\label{2.2}
\|v^\lambda\|_{\mathcal{X}^{-1}} = \int_{|\xi| \leq
1}|\xi|^{-1}|\widehat{v^\lambda}|d\xi + \int_{|\xi| > 1}|\xi|^{-
2}|\widehat{\sqrt{- \Delta} v^\lambda}|d\xi \leq
C\|v^\lambda\|_{H^1}
\end{equation}
and
\begin{eqnarray}\label{2.3}
&&\int_0^{T_\lambda}\|v^\lambda(t)\|_{\mathcal{X}^1}dt =
  \int_0^{T_\lambda}\int_{|\xi| \leq 1}|\xi|
  |\widehat{v^\lambda}|d\xi dt\\\nonumber
&&\quad +\ \int_0^{T_\lambda}\int_{|\xi| > 1} |\xi|^{- 2}
  |\widehat{\sqrt{- \Delta}^3 v^\lambda}|d\xi dt \leq
  C\int_0^{T_\lambda}\|v^\lambda\|_{H^3}ds.
\end{eqnarray}
Hence, one has $v^\lambda \in L^\infty(0, T_\lambda; \mathcal{X}^{-1})
\cap L^1(0, T_\lambda; \mathcal{X}^1)$ (without uniform norms bound in
$\lambda$).
 Next, let us take the
Fourier transform of \eqref{NS} to get
\begin{equation}\label{2.4}
\begin{cases}
\partial_t\widehat{v^\lambda} - i\int \widehat{v^\lambda}(\eta)\otimes
  \widehat{v^\lambda}(\xi - \eta)d\eta\cdot\xi - i\xi \widehat{p^\lambda} + \mu|\xi|^2 \widehat{v^\lambda},\\[-3mm]\\
\xi\cdot \widehat{v^\lambda} = 0.
\end{cases}
\end{equation}
From \eqref{2.4}, \eqref{2.2} and
\eqref{2.3}, we deduce that
\begin{eqnarray}\label{2.5}
&&\partial_t\int|\xi|^{-1}|\widehat{v^\lambda}|d\xi  +
  \mu\int|\xi||\widehat{v^\lambda}|d\xi\\\nonumber
&&= \frac{i}{2}\iint\Big(\big[\widehat{v^\lambda}(\eta)
  \cdot|\widehat{v^\lambda}(\xi)|^{-1}\overline{\widehat{v^\lambda}(\xi)}\big]
  \widehat{v^\lambda}(\xi - \eta)\\\nonumber
&&\quad -\ \big[\overline{\widehat{v^\lambda}}(\eta)
  \cdot|\widehat{v^\lambda}(\xi)|^{-1}\widehat{v^\lambda}(\xi)\big]
  \overline{\widehat{v^\lambda}}(\xi - \eta)\Big)\cdot|\xi|^{-1}\xi d\eta d\xi\\\nonumber
&&\leq \int\int|\widehat{v^\lambda}(\eta)||\widehat{v^\lambda}(\xi
  - \eta)|d\eta d\xi\\\nonumber
&&\leq \frac{1}{2}\int\int\big(|\eta|^{-1}|\xi -
  \eta| + |\eta||\xi - \eta|^{-1}\big)|\widehat{v^\lambda}(\eta)||\widehat{v^\lambda}
  (\xi - \eta)|d\eta d\xi\\\nonumber
&&\leq \int|\xi|^{-1}|\widehat{v^\lambda}(\xi)|d\xi
  \int|\xi||\widehat{v^\lambda}(\xi)|d\xi.
\end{eqnarray}
By \eqref{1.4} and \eqref{2.2}, we see that
$\|v^\lambda(t)\|_{\mathcal{X}^{-1}} < \mu$ at least for a very short
time internal $[0, \delta]$ with $0 < \delta < T_\lambda$.
Consequently, on such a time internal, one has
\begin{eqnarray}\nonumber
\partial_t\|v^\lambda\|_{\mathcal{X}^{-1}} \leq 0,\quad {\rm hence}\
\|v^\lambda\|_{\mathcal{X}^{-1}} \leq \|v_0\|_{\mathcal{X}^{-1}} <
\mu.
\end{eqnarray}
A continuity argument in the time variable yields that
\begin{eqnarray}\nonumber
\|v^\lambda(t)\|_{\mathcal{X}^{-1}} \leq
\|v_0\|_{\mathcal{X}^{-1}} < \mu
\end{eqnarray}
for all $t \in [0, T_\lambda)$. We then apply \eqref{2.5} once more
to derive that
\begin{eqnarray}\label{2.6}
\|v^\lambda(t)\|_{\mathcal{X}^{-1}}  + \big(\mu -
\|v_0\|_{\mathcal{X}^{-1}}\big)\int_0^t\|v^\lambda(s)\|_{\mathcal{X}^1}ds
\leq \|v_0\|_{\mathcal{X}^{-1}},\quad {\rm for\ all}\ t \in [0,
T_\lambda).
\end{eqnarray}

As a bi-product of  \eqref{2.6}, one obtains that
\begin{eqnarray}\nonumber
\int_0^{T_\lambda}\|\nabla v^\lambda(t)\|_{L^\infty}dt \leq
\int_0^{T_\lambda}\|v^\lambda(t)\|_{\mathcal{X}^1}dt \leq
\frac{\|v_0\|_{\mathcal{X}^{-1}}}{\mu -
\|v_0\|_{\mathcal{X}^{-1}}}.
\end{eqnarray}
The standard energy method ( \cite{MB} ) gives that
\begin{eqnarray}\nonumber
\|v^\lambda(t)\|_{H^k} \leq
\|v_0\|_{H^k}\exp\big\{c_k\int_{0}^{T_\lambda}\|\nabla
v^\lambda(s)\|_{L^\infty}ds\big\} \leq
\|v_0\|_{H^k}\exp\big\{\frac{c_k\|v_0\|_{\mathcal{X}^{-1}}}{\mu -
\|v_0\|_{\mathcal{X}^{-1}}}\big\}
\end{eqnarray}
for all $0 \leq t < T_\lambda$ and all $k > 0$. The latter implies that
$T_\lambda = \infty$. Moreover, one has the uniform estimate for
$v^\lambda$:
\begin{equation}\label{2.7}
\sup_{0 \leq t < \infty}\Big(\|v^\lambda(t)\|_{\mathcal{X}^{-1}} +
\big(\mu - \|v_0\|_{\mathcal{X}^{-1}}\big)
\int_0^t\|v^\lambda\|_{\mathcal{X}^1}ds\Big) \leq
\|v_0\|_{\mathcal{X}^{-1}},
\end{equation}
under the condition \eqref{1.4}.

The estimate \eqref{2.7} implies that there exists a subsequence
of $\{v^\lambda\}$ (we will still denote it by $\{v^\lambda\}$)
such that, as $\lambda \rightarrow 0$,
\begin{equation}\label{2.8}
v^\lambda \rightharpoonup v \quad {\rm in}\ L^1(\mathbb{R}_+;
\mathcal{X}^{1}),\quad  v^\lambda \rightharpoonup v \quad {\rm
weakly^\ast\  in}\ L^\infty(\mathbb{R}_+; \mathcal{X}^{- 1})
\end{equation}
for some
\begin{equation}\label{2.9}
v \in L^\infty(\mathbb{R}_+; \mathcal{X}^{-1}) \cap
L^1(\mathbb{R}_+; \mathcal{X}^1).
\end{equation}
For the initial data, we note that
\begin{eqnarray}\nonumber
&&\|v_0^\lambda - v_0\|_{X^{ -1}} = \int_{|\xi| \leq
  M}|\xi|^{- 1}|\widehat{\zeta}(\lambda\xi) -
  1||\widehat{v_0}(\xi)|d\xi\\\nonumber
&&\quad +\ \int_{|\xi| >  M}|\xi|^{- 1}
  |\widehat{\zeta}(\lambda\xi) -
  1||\widehat{v_0}(\xi)|d\xi\\\nonumber
&&\leq 2\sup_{|\eta| \leq \lambda M}
  |\widehat{\zeta}(\eta) - 1|\int_{|\xi| <  M}|\xi|^{- 1}|\widehat{v_0}(\xi)|d\xi +
  2\int_{|\xi| >  M}|\xi|^{- 1}|\widehat{v_0}(\xi)|d\xi.
\end{eqnarray}
By taking $M = \lambda^{- \frac{1}{2}}$, and using the identity
$\widehat{\zeta}(0) = \int\zeta(x)dx = 1$, one concludes that
\begin{eqnarray}\label{2.10}
\|v_0^\lambda - v_0\|_{\mathcal{X}^{ -1}} \rightarrow 0\quad {\rm
as}\ \ \lambda \rightarrow 0.
\end{eqnarray}

To show the strong convergence of $v^\lambda$, we proceed
similarly as in \eqref{2.5} to calculate that
\begin{eqnarray}\nonumber
&&\partial_t\int|\xi|^{-1}|\widehat{v^{\lambda_1}} -
  \widehat{v^{\lambda_2}}|d\xi + \mu\int|\xi|
  |\widehat{v^{\lambda_1}} - \widehat{v^{\lambda_2}}|d\xi\\\nonumber
&&\leq \iint\big(|\widehat{v^{\lambda_1}}(\eta)| +
  |\widehat{v^{\lambda_2}}(\eta)|\big)|\widehat{v^{\lambda_1}}(\xi -
  \eta)  - \widehat{v^{\lambda_2}}(\xi - \eta)|d\eta d\xi\\\nonumber
&&\leq \frac{1}{2}\iint\big(|\eta|^{-1}|\xi -
  \eta| + |\eta||\xi - \eta|^{-1}\big)\\\nonumber
&&\quad \times\ \big(|\widehat{v^{\lambda_1}}(\eta)| +
  |\widehat{v^{\lambda_2}}(\eta)|\big)|\widehat{v^{\lambda_1}}(\xi -
  \eta)  - \widehat{v^{\lambda_2}}(\xi - \eta)|d\eta d\xi\\\nonumber
&&\leq \frac{1}{2}\big(\|v^{\lambda_1}\|_{\mathcal{X}^{-
  1}} + \|v^{\lambda_1}\|_{\mathcal{X}^{- 1}}\big)\|v^{\lambda_1}
  - v^{\lambda_1}\|_{\mathcal{X}^1}\\\nonumber
&&\quad +\ \frac{1}{2}\big(\|v^{\lambda_1}\|_{\mathcal{X}^1} +
  \|v^{\lambda_1}\|_{\mathcal{X}^1}\big)\|v^{\lambda_1}
  - v^{\lambda_1}\|_{\mathcal{X}^{- 1}}.
\end{eqnarray}
Combining the above with \eqref{2.7}, we obtain that
\begin{eqnarray}\nonumber
&&\partial_t\int|\xi|^{-1}|\widehat{v^{\lambda_1}} -
  \widehat{v^{\lambda_2}}|d\xi + \big(\mu - \|v_0\|_{\mathcal{X}^{-1}}\big)\int|\xi|
  |\widehat{v^{\lambda_1}} - \widehat{v^{\lambda_2}}|d\xi\\\nonumber
&&\leq \frac{1}{2}\big(\|v^{\lambda_1}\|_{\mathcal{X}^1} +
  \|v^{\lambda_1}\|_{\mathcal{X}^1}\big)\|v^{\lambda_1}
  - v^{\lambda_1}\|_{\mathcal{X}^{- 1}}.
\end{eqnarray}
The latter implies further that
\begin{equation}\label{2.11}
\begin{cases}
\|v^{\lambda_1}(t) - v^{\lambda_1}(t)\|_{\mathcal{X}^{- 1}} \leq
  \|v_0^{\lambda_1} - v_0^{\lambda_1}\|_{\mathcal{X}^{-
  1}}\exp\big\{\frac{\|v_0\|_{\mathcal{X}^{-1}}}{\mu -
  \|v_0\|_{\mathcal{X}^{-1}}}\big\},\\[-4mm]\\
\big(\mu -
\|v_0\|_{\mathcal{X}^{-1}}\big)\int_0^\infty\|v^{\lambda_1} -
  v^{\lambda_1}\|_{\mathcal{X}^{1}}dt
  \leq \|v_0^{\lambda_1} - v_0^{\lambda_1}\|_{\mathcal{X}^{-
  1}}\exp\big\{\frac{\|v_0\|_{\mathcal{X}^{-1}}}{\mu -
  \|v_0\|_{\mathcal{X}^{-1}}}\big\}.
\end{cases}
\end{equation}
Combining \eqref{2.10} and \eqref{2.11}, we conclude that
$\{v^{\lambda}\}$ is a Cauchy sequence in $L^\infty(\mathbb{R}_+;
\mathcal{X}^{-1}) \cap L^1(\mathbb{R}_+; \mathcal{X}^1)$ and the
convergence in \eqref{2.8} is a strong one. In fact, \eqref{2.11}
also yields the uniqueness of solutions in the space
$L^\infty(\mathbb{R}_+; \mathcal{X}^{- 1}) \cap L^1(\mathbb{R}_+;
\mathcal{X}^1)$ under the assumption \eqref{1.4}.

To get the further time regularity of $v(t, x)$, we come back to
the equation \eqref{2.4}. We claim that $\partial_tv^\lambda$'s are
uniformly bounded in $L^1(\mathbb{R}_+; \mathcal{X}^{- 1})$. Indeed,
by \eqref{2.7}, it is obvious that $\Delta v^\lambda$ are
uniformly bounded in $L^1(\mathbb{R}_+; \mathcal{X}^{- 1})$.
Moreover, as calculated in \eqref{2.5}, one has
\begin{eqnarray}\nonumber
&&\|\nabla\cdot(v^\lambda\otimes v^\lambda)\|_{L^1(\mathcal{X}^{-
  1})} \leq \int_0^\infty\iint|\widehat{v^\lambda}(\eta)|
  |\widehat{v^\lambda}(\xi - \eta)|d\xi d\eta dt\\\nonumber
&&\leq \int_0^t\|v^\lambda\|_{\mathcal{X}^{-
  1}}\|v^\lambda\|_{\mathcal{X}^1}dt \leq \sup_t\|v^\lambda(t)\|_{\mathcal{X}^{-
  1}}\|v^\lambda\|_{L^1(\mathcal{X}^1)}.
\end{eqnarray}
The pressures can be treated by the same way.
We hence proved the claim
\begin{eqnarray}\label{2.12}
\partial_tv^\lambda \in L^1(\mathbb{R}_+; \mathcal{X}^{- 1}),\quad
\|v^\lambda\|_{L^1(\mathcal{X}^1)} \leq C_0\|v_0\|_{\mathcal{X}^{-
1}}\big(1 + \|v_0\|_{\mathcal{X}^{- 1}}\big).
\end{eqnarray}
The latter allows us to improve \eqref{2.9} and to finally conclude
\begin{equation}\label{2.13}
v \in C(\mathbb{R}_+; \mathcal{X}^{-1}) \cap L^1(\mathbb{R}_+;
\mathcal{X}^1),\quad
\partial_tv \in L^1(\mathbb{R}_+; \mathcal{X}^{- 1}).
\end{equation}

\begin{rem}
Let us remark that a similar estimate as \eqref{2.5} also implies
a type of Beale-Kato-Majda's criterion of Navier-Stokes equations:
if
\begin{equation}\nonumber
\int_0^T\int|\widehat{\omega}(\xi)|d\xi dt < \infty,
\end{equation}
then a smooth solution on $[0, T)$ can be extended to $[0, T +
\delta)$ for some $\delta > 0$. In fact,
\begin{eqnarray}\nonumber
\partial_t\int|\xi|^s|\widehat{v}|d\xi  +
  \mu\int|\xi|^{s + 2}|\widehat{v}|d\xi \leq \iint|\xi|^{s + 1}|
  \widehat{v}(\eta)||\widehat{v}(\xi
  - \eta)|d\eta d\xi.
\end{eqnarray}
Noting $v = \Delta^{-1}{\rm curl}\omega$ and taking $s = - 1$ and
0, one gets that
$$\|v(t)\|_{\mathcal{X}^{-1}} \leq \|v_0\|_{\mathcal{X}^{-1}}
\exp\big\{\int_0^t\int|\widehat{\omega}(\xi)|d\xi ds\big\} <
\infty$$ and
$$\|v(t)\|_{\mathcal{X}^0} \leq \|v_0\|_{\mathcal{X}^0}\exp\big\{
2\int_0^t\int|\widehat{\omega}(\xi)|d\xi ds\big\} < \infty$$ for
all $0 \leq t < T$. Then for $s > 0$ one has
\begin{eqnarray}\nonumber
&&\partial_t\int|\xi|^s|\widehat{v}|d\xi  + \mu\int|\xi|^{s + 2}
  |\widehat{v}|d\xi \leq \iint_{|\xi| \leq M_s}|\xi|^{s + 1}|
  \widehat{v}(\eta)||\widehat{v}(\xi
  - \eta)|d\eta d\xi\\\nonumber
&&\quad +\ \epsilon_s\iint_{|\xi| > M_s}(1 + |\eta|^{s + 2} + |\xi
  - \eta|^{s + 2})\widehat{v}(\eta)||\widehat{v}(\xi
  - \eta)|d\eta d\xi\\\nonumber
&&\leq (M_s^{s + 1} + \epsilon_s)\|v\|_{\mathcal{X}^{-1}}\|v
  \|_{\mathcal{X}^1} + 2\epsilon_s\int|\widehat{v}|d\xi\int
  |\xi|^{s + 2}|\widehat{v}(\xi)|d\xi.
\end{eqnarray}
Here $\epsilon_s > 0$ is a small constant such that
$2\epsilon_s\int|\widehat{v_0}|d\xi\exp\{2\int_0^T\int|\widehat{\omega}(\xi)|d\xi
dt\} < \mu$ which implies $2\epsilon_s\int|\widehat{v}|d\xi <
\mu$. Then $M_s > 1$ is chosen to be a large constant such that
$|\xi|^{s + 1} \leq \epsilon_s(1 + |\eta|^{s + 2} + |\xi -
\eta|^{s + 2})$ for $|\xi| > M_s$.
 Consequently, one has $\|v(t)\|_{C^s}  < \infty$. We
emphasis that the bounds of $\|v\|_{C^s}$ are only one exponential
in terms of $\int_0^T\int|\widehat{\omega}(\xi)|d\xi dt$.
\end{rem}

\section*{Acknowledgement}
The authors would like to thank the anonymous referee for
suggesting us the counterexample at the end of the introduction.
Zhen Lei was in part supported by NSFC (grants No. 10801029 and
10911120384), FANEDD, Shanghai Rising Star Program (10QA1400300),
SGST 09DZ2272900 and SRF for ROCS, SEM. Fanghua Lin is partially
supported by an NSF grant, DMS0700517. Part of the work was
carried out while Zhen Lei was visiting the Courant Institute.


\end{document}